\documentclass[11pt,reqno,a4paper]{amsart}

\usepackage{amsmath,amssymb,amsthm,amsaddr}
\usepackage{mathtools}
\usepackage{enumitem}
\usepackage[margin=3cm,top=4cm,bottom=4cm,footskip=1cm,headsep=2cm]{geometry}
\usepackage[utf8]{inputenc}

\usepackage{amsthm}
\usepackage[british]{babel}
\usepackage{xcolor}

\def\div{\operatorname{div}}
\def\dt{\partial_t}
\def\TT{\mathbb{T}}
\def\RR{\mathbb{R}}
\def\dx{\partial_x}

\def\eps{\varepsilon}
\def\H{\mathcal{H}}

\newtheorem{lemma}{Lemma}

\newtheorem{theorem}[lemma]{Theorem}

\theoremstyle{definition}
\newtheorem{remark}[lemma]{Remark}

\DeclarePairedDelimiter{\norm}{\|}{\|}

\def\softd{{\leavevmode\setbox1=\hbox{d}%
		\hbox to 1.05\wd1{d\kern-0.25ex{\char039}\hss}}}%cstocs

\title[Parameter estimation in the Cahn-Hilliard equation]{On uniqueness and stable estimation of multiple parameters in the Cahn-Hilliard equation}

\author{Aaron Brunk$^1$, Herbert Egger$^{2,3}$, Oliver Habrich$^3$}

\address{$^1$Institute of Mathematics, Johannes Gutenberg-University Mainz, Germany}
\address{$^2$Johann Radon Institute for Computational and Applied Mathematics, Linz, Austria}
\address{$^3$Institute of Numerical Mathematics, Johannes Kepler University Linz, Austria}

\begin{document}

\begin{abstract}
We consider the identifiability and stable numerical estimation of multiple parameters in a Cahn-Hilliard model for phase separation.  
Spatially resolved measurements of the phase fraction are assumed to be accessible, with which the identifiability of single and multiple parameters up to certain scaling invariances is established. 
A regularized equation error approach is proposed for the stable numerical solution of the parameter identification problems, and convergence of the regularized approximations is proven under reasonable assumptions on the data noise. 
The viability of the theoretical results and the proposed methods is demonstrated in numerical tests.
\end{abstract}

\maketitle

\begin{quote}

\medskip

\end{quote}

\bigskip

\section{Introduction}

The modelling and simulation of phase transformation processes is of interest in many applications, e.g., in the spinodal decomposition of binary alloys and fluid mixtures \cite{elliott1989cahnhilliard,tanaka1998spinodal}, or in the modelling of tumour growth \cite{cristini2010multiscale,garcke2016cahnhilliard}. 
One of the basic mathematical models arising in such applications is the Cahn-Hilliard system 
\begin{align}
    \dt \phi  &=  \div\left( b(\phi)\nabla \mu \right), \label{eq:ch1}  \\
    \mu &   = - \gamma \Delta \phi + f(\phi). \label{eq:ch2}
\end{align}
Here $\phi$ is the phase fraction of the mixture, $\mu$ the chemical potential, $b(\phi)$ a concentration dependent mobility, $\gamma$ an interface parameter, and $f(\phi)=F'(\phi)$ the derivative of a double well potential $F(\phi)$ whose minima characterize the favourable mixing ratios. 

Together with appropriate boundary conditions, the Cahn-Hilliard equation describes the gradient flow of a conserved quantity, i.e., the free energy 
\begin{align} \label{eq:ch3}
    \mathcal{E}(\phi) = \int_\Omega \frac{\gamma}{2} |\nabla \phi|^2 + F(\phi) \, dx 
\end{align}
is dissipated over time while the total amount $\mathcal{M}(\phi)=\int_\Omega \phi \, dx$ of substance is preserved.
These are key ingredients for establishing existence of solutions and they also guarantee the thermodynamic consistency of the model \cite{BarrettBlowey99,Hohenberg}. 
The Cahn-Hilliard equation, however, is still a phenomenological model describing the macroscopic behaviour of the system under consideration, and a careful calibration of the model parameters is required to obtain quantitative agreement with more detailed microscopic descriptions or experimental data \cite{cristini2010multiscale,hulikal2016experimental}.
Finding the model parameters in  \eqref{eq:ch1}--\eqref{eq:ch2} from observations of the solution amounts to a typical parameter estimation problem in a nonlinear system of partial differential equations; see \cite{banks1989book,isakov2017book} for an introduction and references. 

Parameter identification problems in nonlinear parabolic equations and methods for their stable solution have been studied intensively in the literature, in particular, in the context of heat transfer and porous medium flow. 
In \cite{cannon1980inverse,chavent1974identification,duemmel1988uniqueness}, the identification of the nonlinear conductivity function $a(u)$ in $\dt u = \div(a(u) \nabla u)$ has been addressed. Uniqueness results in one and multiple space dimensions have been derived, and an output least squares method has been used for the numerical solution. 
Equation error methods are proposed and thoroughly analyzed in \cite{cao2006natural,hanke1999equation} for the stable numerical solution.
These approaches have been developed for related linear elliptic problems in \cite{acar1993identification,alessandrini1986identification,kaerkkaeinen1997equation,kohn1988variational}; also see \cite{aljamal2012stability,kaltenbacher2002saddle}.
The simultaneous identification of multiple parameters in nonlinear elliptic and parabolic problems has been investigated, for instance, in \cite{bitterlich2002efficient,bukshtynov2011optimal,egger2014simultaneous,fatullayev2007numerical,ou2008inverse}, and uniqueness questions as well as numerical methods for the stable solution have been studied. 
Related results have also been derived in the context of chemotaxis \cite{egger2015identification,fister2008identification}. 
The recent work \cite{kahle2020parameter,kahle2019bayesian} addresses the identification of multiple scalar parameters in a phase-field model for tumour growth, which is an extended version of the Cahn-Hilliard system above.

In this paper, we study, theoretically and numerically, the identification of 
\begin{align} \label{eq:ch4}
\gamma, \quad b(\cdot), \quad f(\cdot)
\end{align}
in the nonlinear parabolic system \eqref{eq:ch1}--\eqref{eq:ch2} from distributed measurements of $\phi$. Note that such observations of the phase fraction $\phi$ are available from simulations of microscopic models or experimental investigations \cite{hulikal2016experimental}, while similar information about the chemical potential $\mu$ is typically not directly accessible in practice. 
One may therefore classify the parameter identification problem as one with \emph{incomplete data}.

\subsection*{Outline and main results}
We start with introducing our notation and some preliminary results in Section~\ref{sec:prelim}, and identify invariances of the problem with respect to certain scalings, which characterize an inherent non-uniqueness of the parameter identification problem. This allows us to eliminate the parameter $\gamma$ from our consideration. 
In Section~\ref{sec:pid}, we then investigate the independent identification of the potential $f(\cdot)$ and the mobility $b(\cdot)$, as well as the simultaneous identification of both parameter functions, and we establish uniqueness results for the corresponding inverse problems under certain observability conditions; see \cite{alessandrini1986identification,richter1981numerical} for similar conditions in the context of linear parabolic and elliptic equations.
In Section~\ref{sec:reg}, we turn to the stable numerical solution, for which we employ an equation error approach in the spirit of \cite{egger2015identification,hanke1999equation}. This reduces the parameter identification problems to linear ill-posed operator equations with perturbed operators, whose stable solution is accomplished by standard regularization methods; related analysis can be found in \cite{cao2006natural,egger2015identification}.
For illustration of our results, some numerical tests are presented in Section~\ref{sec:num}, and the presentation closes with a short discussion.

\section{preliminaries} \label{sec:prelim}

Let us start with introducing our basic assumptions on the computational domain and the coefficients. For ease of presentation, we consider \eqref{eq:ch1}--\eqref{eq:ch2} on a $d$-dimensional cube and complemented by periodic boundary conditions, i.e., 
\begin{itemize}\itemsep1ex
\item[(A0)] $\Omega \simeq \TT^d$, is the $d$-dimensional torus for $d=1,2,3$. \\
Moreover, functions defined on $\Omega$ are always assumed to be periodic.
\end{itemize}
We further impose the following assumptions on the model parameters
\begin{itemize}\itemsep1ex
\item[(A1)] $\gamma>0$ is a positive constant;
\item[(A2)] $b: \RR \to \RR_+$ satisfies $b\in C^2(\mathbb{R})$ with $0< b_1\leq b(s) \leq b_2$, $\norm{b'}_\infty\leq b_3$, $\norm{b''}_\infty\leq b_4$;
\item[(A3)] $f(s) = F'(s)$ with $F\in C^4(\mathbb{R})$ such that $F(s),F''(s)\geq -f_1$, for some $f_1\geq0$. Furthermore, $F$ and its derivatives are bounded by $|F^{(k)}(s)| \le f_2^{(k)} + f_3^{(k)} |s|^{4-k}$ for $0 \le k \le 4$ with constants $f_2^{(k)}, f_3^{(k)}\geq 0$.
\end{itemize}
These standard assumptions allow to establish existence, uniqueness and regularity for solutions of the Cahn-Hilliard system; see Lemma~\ref{lem:weak} below.
By $L^p(D)$ and $H^{k}(D)$, we denote the standard Lebesgue and Sobolev spaces over some manifold  $D$.
For $D=\Omega \simeq \TT^d$, the functions in these spaces are assumed to be periodic. 
By $L^p(X)=L^p(0,T;X)$, we denote the Bochner spaces of $L^p$ integrable functions $g:[0,T] \to X$ with values in some Banach space~$X$. 
All spaces are equipped with their standard norms; see e.g. \cite{Evans}. 

\subsection{Cahn-Hilliard equation}
For later reference, let us recall some well-known results about existence, uniqueness and regularity of solution to the Cahn-Hilliard equation.

\begin{lemma} \label{lem:weak}
Let (A0)--(A3) hold. Then for any $\phi_0 \in H^1(\Omega)$,
there exists at least one periodic weak solution 
\begin{align*}
(\phi,\mu)\in H^1(0,T;H^{-1}(\Omega))\cap L^2(0,T;H^3(\Omega))\times L^2(0,T;H^1(\Omega))
\end{align*}
of problem \eqref{eq:ch1}--\eqref{eq:ch2} with initial value $\phi(0)=\phi_0$.  
If $\phi_0 \in H^3(\Omega)$, and $T$ sufficiently small in dimension $d=3$, then
\begin{align*}
    \|\phi\|_{L^\infty(H^3)} + \|\dt \phi\|_{L^\infty(H^{-1})} \le C_{T},
\end{align*}
with constant $C_{T}$ depending only on the bounds for the coefficient, the domain $\Omega$, the time horizon $T$, and the bounds for the initial value.
Moreover, the weak solution is unique.
\end{lemma}
A proof of these assertions is obtained by standard energy methods and Galerkin approximations; see \cite{Boyer,ElliotGarcke2000,Huang2009} for details and related results. With similar arguments and further smoothness assumptions on the parameter functions and initial value, also higher regularity of the solutions can be established.

\begin{remark} \label{rem:stationary}
If $\phi_0 \in H^3(\Omega)$ and, in addition, the function $F$ is analytic, then the phase fraction $\phi$ converges to an equilibrium distribution $\phi_\infty$, 
more precisely
\begin{align*}
\|\phi(\cdot,t)-\phi_\infty\|_{H^3(\Omega)} 
\leq C_\infty (1+t)^{-\frac{\theta}{1-2\theta}}, \forall t\geq 1,
\end{align*}
with a constant $C_\infty$ depending on the parameters, the initial data, on $\phi_\infty$, and the rate constant $\theta\in(0,\frac{1}{2})$  depending on $\phi_\infty$;
we refer to \cite{Rybka,Huang2009} for details and proofs.
\end{remark}

\subsection{Scaling invariances}

Before we turn to a detailed statement and analysis of the parameter identification problems, let us highlight the following canonical invariances, which characterize the inherent non-uniqueness for the parameter identification problems.

\begin{lemma} \label{lem:invariance}
Let $(\phi,\mu)$ be a periodic solution of  \eqref{eq:ch1}--\eqref{eq:ch2} for parameters $(\gamma,b,f)$. Then for any $c \in \RR$ and $d>0$, the tuple $(\hat \phi,\hat \mu)=(\phi,\mu/d+c)$ is a solution of \eqref{eq:ch1}--\eqref{eq:ch2} for parameters 
\begin{align*}
\hat \gamma =\gamma/d, 
\qquad 
\hat b(s) = d \cdot  b(s), 
\qquad \text{and} \qquad 
\hat f(s) = f(s) / d + c.
\end{align*}
\end{lemma}
\begin{proof}
It is easy to see that a constant scaling $b \to b \cdot d$ and $\mu \to \mu/d$ does not perturb the validity of \eqref{eq:ch1}, and rescaling $\gamma \to \gamma/d$, $f \to f/d$ restores the validity of \eqref{eq:ch2}. 
Further, note that a constant shift $\mu \to \mu+c$ and $f \to f+c$ leaves the equations \eqref{eq:ch1}--\eqref{eq:ch2} valid. 
A combination of the two scalings already yields the result.
\end{proof}

\begin{remark} \label{rem:invariance}
Using distributed observations of the phase fraction $\phi$ only, the parameters $\gamma$, $b(\cdot)$, $f(\cdot)$ can be identified at most up to the above invariant scalings. 
Without loss of generality, we therefore assume $\gamma>0$ to be given in the following, and we consider the identification of the functions $b(\cdot)$ and $f(\cdot)=F'(\cdot)$, the latter up to a constant. 
\end{remark}

\section{Identifiability results} \label{sec:pid}

In the following, we first study the separate identification of $f(\cdot)$ and $b(\cdot)$, and then turn to the simultaneous identification of both parameters. 
According to Remark~\ref{rem:invariance}, we assume that $\gamma>0$ is known, and then expect that $b(\cdot)$ can be identified uniquely from distributed measurements of $\phi$, while $f(\cdot)$ can be determined up to a constant shift.

\subsection{Identification of $f(\cdot)$}

We eliminate $\mu$ by inserting \eqref{eq:ch2} into \eqref{eq:ch1}. The resulting equation is then multiplied by a periodic test function $v$ and integrated over the domain~$\Omega$. After integration-by-parts and using the periodicity of the $\phi$ and $v$, we see that 
\begin{align} \label{eq:f1}
\int_\Omega b(\phi) f'(\phi) \nabla \phi \cdot \nabla v \, dx &=  \int_\Omega b(\phi)\gamma \nabla \Delta\phi \cdot \nabla v \, dx -\int_\Omega \dt \phi \, v \, dx.
\end{align}
Since we assumed $b(\cdot)$ strictly positive, we may define $v = \int_0^\phi w(s)/b(s) \, ds$ for any smooth periodic test function $w$, and the above variational identity leads to
\begin{align} \label{eq:f2}
\int_\Omega f'(\phi) w(\phi) |\nabla \phi|^2 \, dx &=  \int_\Omega  w(\phi) \gamma \nabla \Delta \phi \cdot \nabla \phi \, dx -\int_\Omega  \int_0^\phi w(s)/b(s) \, ds \, \dt \phi \, dx.  
\end{align}
From this identity, we can immediately deduce the following result.
\begin{theorem} \label{thm:uniq-f}
Let (A0)--(A3) hold and $\phi$ be a smooth solution of \eqref{eq:ch1}--\eqref{eq:ch2} on $[0,T] \times \Omega$. 
Further assume that $\gamma>0$ and $b(\cdot)$ are known.
Then $f'(\cdot)$ is uniquely determined on $R_t := \{s = \phi(x,t) : x \in \Omega\}$ from observations of $\phi(\cdot,t)$ and $\dt \phi(\cdot,t)$ on $\Omega$ for $0 \le t \le T$.
\end{theorem}
\begin{proof}
Without mentioning explicitly, we always consider a specific time point $t \in [0,T]$ in the following. 
Assume that $f_1$, $f_2$ are two functions leading to the same solution $\phi$. Then 
\begin{align} \label{eq:f3}
    \int_\Omega (f_1'(\phi) - f_2'(\phi)) \,  w(\phi) \, |\nabla \phi|^2 \, dx &= 0.
\end{align}
We define $W(s) = f_1(s) - f_2(s)$ and $w(s) = W'(s)$, and note that $\nabla W(\phi) = w(\phi) \nabla \phi$. 
Using this $w(\cdot)$ as a test function in the above variational identity, we obtain
\begin{align*}
\int_\Omega    |\nabla (f_1(\phi) - f_2(\phi))|^2 dx 
= \int_\Omega |f_1'(\phi) - f_2'(\phi)|^2 |\nabla \phi|^2 dx 
= 0.
\end{align*}
This implies that $f_1(\phi) - f_2(\phi)$ is constant on $\Omega$, which already yields the claim.  
\end{proof}

\begin{remark}
If we have measurements of $\phi$ on a space-time cylinder $\Omega \times [t_1,t_2]$ with $t_1 < t_2$, then we also know $\dt \phi$ on this set, and can determine $f'(\cdot)$ on $R_{[t_1,t_2]}=\bigcup_{t_1 \le t \le t_2} R_t$.  
\end{remark}

Let us note that the last term in equation \eqref{eq:f2} depends on $b(\cdot)$, and we therefore had to assume that the mobility function is known in order to derive \eqref{eq:f3}. 
Alternatively, the dependence on $b(\cdot)$ in equation~\eqref{eq:f2} disappears, if we assume that $\dt \phi=0$ on $\Omega$, which corresponds to an equilibrium situation; see Remark~\ref{rem:stationary}.
In this case, $f(\cdot)$ can be determined up to constants without knowledge of the mobility. 
This leads to the following result.

\begin{theorem}\label{thm:ident-f-equi}
Let 
(A0)--(A3) hold and
$(\phi_\infty,\mu_\infty)$ be an equilibrium for \eqref{eq:ch1}--\eqref{eq:ch2}, i.e., 
\begin{align}
    0 &= \div(b(\phi_\infty) \nabla \mu_\infty) \label{eq:eq1}\\
    \mu_\infty &= -\gamma \Delta \phi_\infty +
     f(\phi_\infty). \label{eq:eq2}
\end{align}
Further assume that $\gamma>0$ to be known.
Then the function $f'(\cdot)$ is determined uniquely on $R_\infty=\{s =\phi_\infty(x) : x \in \Omega\}$ from knowledge of $\phi_\infty$. 
\end{theorem}
\begin{proof}
In principle, the assertion follows directly from the previous theorem, but we present a more direct proof here.   
From \eqref{eq:eq1} we conclude that 
\begin{align*}
0 
&= \int_\Omega \div(b(\phi_\infty) \nabla \mu_\infty) \mu_\infty \, dx 
 = -\int_\Omega b(\phi_\infty) |\nabla \mu_\infty|^2 dx. 
\end{align*}
Since the mobility $b(\cdot)$ was assumed to strictly positive, this implies $\mu_\infty \equiv C$ constant. Inserting this into equation \eqref{eq:eq2} then leads to 
\begin{align*}
    f(\phi_\infty) = \gamma \Delta \phi_\infty - C,
\end{align*}
which already yields the assertion of the lemma.
Let us note that knowledge of the mobility function $b(\cdot)$ was not required here.
\end{proof}

\subsection{Identification of $b(\cdot)$}

We now study the identification of the mobility function $b(\cdot)$, while assuming that the other parameters $\gamma$ and $f(\cdot)$ are known. 
In this case, the chemical potential 
$\mu = -\gamma \Delta \phi + f(\phi)$
is fully determined from observations of $\phi$ already.
We can then rewrite equation \eqref{eq:ch1} as 
\begin{align} \label{eq:ch1b}
    \div (b(\phi) \nabla \mu) = \dt \phi \qquad \text{on } \Omega \times (0,T).
\end{align}
Since $\phi$, $\dt \phi$ and $\mu$ are known at this point, this can be interpreted as a linear operator equation for determining the mobility function $b(\cdot)$. 
Using a similar argument as employed in reference~\cite{egger2015identification}, we now obtain the following conditional identifiability result. 
\begin{theorem}\label{thm:identb}
Let (A0)--(A3) hold and $(\phi,\mu)$ be a smooth solution of \eqref{eq:ch1}--\eqref{eq:ch2} on $[0,T] \times \Omega$. 
Further, assume that $\gamma>0$ and $f(\cdot)$ are known. 
Then $b(\cdot)$ can be determined uniquely from observations of $\phi(\cdot,t)$ and $\dt \phi(\cdot,t)$ on the set
\begin{align*}
\widetilde R_t = \{s = \phi(x,t) : x \in \Omega \ \text{and} \ \nabla \mu(x,t) \ne 0\} \subset R_t.
\end{align*}
\end{theorem}
\begin{proof}
Without further noticing, we always consider a fixed time point $t$ in the following. 
Assume that $\phi$ is the solution of \eqref{eq:ch1}--\eqref{eq:ch2} for the same $\gamma$ and $f(\cdot)$, but for two different mobility functions $b_1(\cdot)$, $b_2(\cdot)$. 
We define $B_+(s) = \max(b_1(s) - b_2(s),0)$, and observe that
\begin{align} \label{eq:ch1c}
    \div ( B_+(\phi) \nabla \mu) = 0 \qquad \text{on }\Omega.
\end{align}
To see this, first, note that the function $B_+(\phi) \nabla \mu$ is weakly differentiable. 
The validity of \eqref{eq:ch1c} on the set $\Omega_+(t) := \{x \in \Omega : B_+(\phi(x,t))>0\} = \{x \in \Omega : b_1(\phi(x,t)) > b_2(\phi(x,t))\}$ then follows by subtracting equation~\eqref{eq:ch1} for $b=b_1$ and $b=b_2$, while validity on $\Omega \setminus \Omega_+$ is trivial, since $B_+(\phi) \nabla \mu \equiv 0$ there. 
We can now multiply equation~\eqref{eq:ch1c} by $\mu$ and integrate over the domain to see that 
\begin{align*}
0 &= \int_\Omega \div (B_+(\phi) \nabla \mu) \mu \, dx 
= -\int_\Omega B_+(\phi) |\nabla \mu|^2 \, dx.
\end{align*}
In the second step, we again employed integration-by-parts and the periodicity of the solutions.
This shows that $B_+(\phi) \equiv 0$ on the set
$\{s=\phi(x,t) \in \tilde R_t : x \in \Omega_+\}$.
In the very same manner, one can verify that $B_-(\phi) = \min(b_1(\phi) - b_2(\phi) , 0)$ vanishes on the remaining set $\{s=\phi(x,t) \in \tilde R_t : x \in \Omega_-\}$,
%$\widetilde R_t \cap \Omega_-$,
where $\Omega_-=\{x \in \Omega : b_1(\phi(x,t))<b_2(\phi(x,t))\}$.
In summary, we thus obtain  $b_1(s)=b_2(s)$ on $\widetilde R_t$.
\end{proof}

\begin{remark}
If data $\phi(x,t)$ are available on a space-time cylinder $\Omega \times [t_1,t_2]$ for $t_1 < t_2$, we also know $\dt \phi$ on this set and can determine $b(\cdot)$ on $\widetilde R_{[t_1,t_2]} = \bigcup_{t_1 \le t \le t_2} \widetilde R_t$. 
Since the function $b(\cdot)$ was assumed smooth, we can take the closure of the sets. 
On the other hand, if $\nabla \mu  \equiv 0$ on $\{(x,t) : \phi(x,t) \in (s_1,s_2)\}$, then $b(s)$ for $s_1 < s < s_2$ obviously has no influence on the evolution of $\phi$ and hence cannot be determined from observations of $\phi$. 
The non-vanishing condition on $\nabla \mu$ contained in the definition of the sets $\widetilde R_t$ and $\widetilde R_{[t_1,t_2]}$ can be understood as an observability condition; compare with \cite{alessandrini1986identification,richter1981numerical}.
\end{remark}

\subsection{Simultaneous identification of $f(\cdot)$ and $b(\cdot)$.}

We will now demonstrate that, under a suitable observability condition, also the simultaneous identification of both parameter functions is possible; see \cite{alessandrini1986identification,richter1981numerical} for similar arguments.
Let $\phi$ be a smooth periodic solution of \eqref{eq:ch1}--\eqref{eq:ch2} and $c(s):= b(s)f'(s)$. Then, from \eqref{eq:f1}, we may conclude that 
\begin{align} \label{eq:bcw}
    \int_\Omega  -b(\phi )\gamma \nabla \Delta \phi \nabla v \, dx + \int_\Omega c(\phi)\nabla \phi \nabla v \, dx     = \int_\Omega \dt \phi \, v \, dx,
\end{align}
for all periodic test functions $v$ and all-time points $t$ under consideration. 
If we assume that $\phi$, $\dt \phi$, and $\gamma>0 $ are known, this can be interpreted as  the variational form of a linear operator equation for determining the two parameter functions $b(\cdot)$ and $c(\cdot)$. 
As a next step, let us define 
\begin{align*}
H(s) = \begin{cases}
0, & s \le 0, \\
1, & s >0, 
\end{cases}
\qquad \text{and} \qquad 
H_\eps(s) = \begin{cases} 
s/\eps, & 0 < s < \eps, \\
H(s), & \text{else},
    \end{cases}
\end{align*}
i.e., the Heaviside function $H(s)$ and its regularized piecewise linear approximations. 
Then by testing \eqref{eq:bcw} with $v=H_\eps(\phi - s)$ and taking the limit $\eps \to 0$, we obtain  
\begin{align} \label{eq:caf}
-b(s) \int_{\{\phi = s\}} \gamma \nabla \Delta \phi \frac{\nabla \phi}{|\nabla \phi|} d\H^{d-1} 
+ c(s) \int_{\{\phi = s\}} |\nabla \phi| \,  d\H^{d-1}
    = \int_{\Omega} \dt \phi \, H(\phi - s) \, dx,
\end{align}
where $d\H^{d-1}$ denotes the $(d-1)$-dimensional Hausdorff measure. 
The integrals on the left-hand side are the same as those appearing in the co-area formula, and thus well-defined for a.e. $s$; see \cite[Sec.~3.4]{evans1991measure} for details.
For every time $t$, this yields a linear equation for the two scalar values $b(s)$ and $c(s)$. This leads to the following result.
\begin{theorem} \label{thm:bc}
Let (A0)--(A3) hold, $\phi$ be a smooth solution of \eqref{eq:ch1}--\eqref{eq:ch2} and $\gamma>0$ known.
Define 
$A_b(s,t) = -\int_{\{\phi(\cdot,t)=s\}} \gamma \nabla \Delta \phi(x,t) \frac{\nabla \phi(x,t)}{|\nabla \phi(x,t)|} \, dx$,
$A_c(s,t) = \int_{\{\phi(\cdot,t)=s\}} |\nabla \phi(x,t)| \, dx$, 
and 
$A(s,t) = \int_\Omega \dt \phi(x,t) H(\phi(x,t) - s) \, dx$,
and assume that $\{(A_b(s,t_i),A_c(s,t_i)): i=1,2\}$ are linearly independent. 
Then $b(s)$, $c(s)$, and $f'(s) = c(s)/b(s)$ are uniquely determined. 
\end{theorem}
\begin{proof}
From equation \eqref{eq:caf} and definition of $A_b$, $A_c$ and $A$, we see that
\begin{align*}
    A_b(s,t_1) \, b(s) + A_c(s,t_1) \, c(s) &= A(s,t_1), \\
    A_b(s,t_2) \, b(s) + A_c(s,t_2) \, c(s) &= A(s,t_2).
\end{align*}
This is a system of two linear equations for determining the two scalar values $b(s)$, $c(s)$, and by assumption, the two equations are linearly independent. 
\end{proof}

\begin{remark}
Theorem~\ref{thm:bc} is a conditional identifiability result.
The required linear independence of the coefficients $A_b$, $A_c$ can, in principle, be checked explicitly using the data. 
If this observability condition is valid, then the two scalar values $b(s)$, $f'(s)$ can be determined uniquely by the linear system \eqref{eq:caf} and then also depend stably on $\phi$ and $\dt \phi$.
\end{remark}

\section{Regularized inversion by equation error methods} 
\label{sec:reg}

We now discuss the stable identification of $f(\cdot)$ and $b(\cdot)$ from observations of $\phi$ also from a numerical point of view. 
We assume that spatially resolved measurements $\phi^\delta$ of $\phi$ are available for particular time steps $t \in [0,T]$, and that
\begin{align}\label{eq:ass-noise}
    \| \phi(\cdot,t) - \phi^\delta(\cdot,t) \|_{H^3(\Omega)} \leq \delta, \qquad \| \dt \phi(\cdot,t) - \dt \phi^\delta(\cdot,t) \|_{H^{-1}(\Omega)} \leq \delta,
\end{align}
with known noise level $\delta$.
In view of the regularity results for the true solution, stated in Lemma~\ref{lem:weak}, such an assumption is realistic after appropriate pre-smoothing of the data~\cite{kaltenbacher2002saddle}. 

\subsection{Equation error approach}

For the numerical solution of the three parameter identification problems which were analyzed in the previous section, we consider an equation error approach~\cite{banks1989book,hanke1999equation}, where the data are directly inserted into the partial differential equations.
This allows to reduce the non-linear parameter identification problems to linear operator equations of the form
\begin{align} \label{eq:linop}
    T^\delta x = y^\delta,
\end{align}
with perturbed operators $T^\delta$ and data $y^\delta$. 
Stable approximations for the solution $x$, i.e., the unknown parameter functions, can then be derived using Tikhonov regularization.
We refer to \cite{cao2006natural,egger2015identification,hanke1999equation} for analysis and examples. 
Before returning to the identification of the parameter functions $f(\cdot)$ and $b(\cdot)$ in \eqref{eq:ch1}--\eqref{eq:ch2}, let us recall the following abstract result, whose proof can be found in \cite{cao2006natural,egger2015identification}, which will serve as a theoretical backup.

\begin{lemma} \label{lem:convergenceapproximations}
Let $T,T^\delta:X \to Y$ be bounded linear operators between Hilbert spaces $X$ and~$Y$. Further let $y\in R(T)$, $y^\delta \in Y$,
and assume that 
\begin{align} \label{eq:condition}
\| T x^\dag - T^\delta x^\dag \|_Y \leq C \delta
\qquad \text{and} \qquad
\|y - y^\delta \| \leq C' \delta 
\end{align}
Then for $\alpha \to 0$ and $\delta^2/\alpha \to 0$, 
the regularized solutions $x_\alpha^\delta$, determined by 
\begin{align}\label{eq:leastsquaresprob}
  \| T^\delta x - y^\delta \|_{Y}^2 + \alpha \| x \|_{X}^2 \to \underset{x\in X}{\min},
\end{align}
converge to the minimum-norm solution $x^\dag$ of $Tx = y$ with $\delta \to 0$.
\end{lemma}
In the following, we show how to transform the parameter identification problems of the previous section into linear operator equations of the form \eqref{eq:linop}. 
Using assumption~\eqref{eq:ass-noise}, we will show that the condition~\eqref{eq:condition} of the previous lemma is satisfied. 
Convergence of the regularized solutions is then guaranteed by the abstract theoretical result above.
As outlined in Remark~\ref{rem:invariance}, the interface parameter $\gamma>0$ is always assumed to be known.

\subsection{Identification of $f(\cdot)$} \label{sec:identify_f}
We first assume $\gamma>0$ and $b(\cdot)$ to be known and study the identification of $f(\cdot)$ in \eqref{eq:ch1}--\eqref{eq:ch2}.
Similar as in Section~\ref{sec:pid}, we eliminate the chemical potential $\mu$ by inserting~\eqref{eq:ch2} into \eqref{eq:ch1}, and define $c(s):=b(s)f'(s)$. This leads to
\begin{align} \label{eq:invf1}
    \div ( c(\phi)\nabla \phi  ) = \gamma \div ( b(\phi) \nabla \Delta \phi) +\dt \phi.
\end{align}
Since the mobility function $b(\cdot)$ is assumed to be known and strictly positive, we can uniquely and stably determine $f'(\cdot)$ from knowledge of $c(\cdot)$.
Like before, we first consider a single time instance $t \in (0,T)$, and we write $\phi$ for $\phi(\cdot,t)$. 
In view of Theorem~\ref{thm:uniq-f}, the parameter function $f'(\cdot)$ can be determined uniquely on the range $R_t$
of values of $\phi$ attained by the data.
A similar statement holds for observations on a whole time interval.

Replacing the exact solution $\phi$ in \eqref{eq:invf1} by the perturbed data $\phi^\delta$ leads to a linear operator equation of the form \eqref{eq:linop}, which can be used to identify the parameter function $c(\cdot)$.
For ease of presentation, we assume $\phi^\delta(\cdot,t) \in (-1,1)$ in the following, which is also satisfied in our numerical tests.
We then define the perturbed operator
\begin{align} \label{eq:Tf}
    T^\delta: H^2(-1,1) \to H^{-1}(\Omega), \quad c(\cdot) \to \div( c(\phi^\delta ) \nabla \phi^\delta ).
\end{align}
The perturbed right hand side is given by $y^\delta := \dt \phi^\delta + \gamma \div ( b(\phi^\delta) \nabla \Delta \phi^\delta )$.
It is not difficult to verify that the operator $T^\delta$ is linear. 
Using integration-by-parts, we see that 
\begin{align*}
  \|T^\delta c\|_{H^{-1}(\Omega)} 
  &= \underset{v\in H^1(\Omega)}{\sup} \frac{(c(\phi^\delta) \nabla \phi^\delta, \nabla v )_{L^2(\Omega)}}{\|v\|_{H^1(\Omega)}}  \\
  &\leq \|c(\phi^\delta) \nabla \phi^\delta \|_{L^2(\Omega)}
  \le \|c\|_{L^\infty(-1,1)} \|\phi^\delta\|_{H^1(\Omega)}.
\end{align*}
From equation~\eqref{eq:ass-noise}, we know that the last term is bounded, and by Sobolev's embedding theorem, we conclude that $\|T^\delta c\|_{H^{-1}(\Omega)} \le C \|c\|_{H^2(-1,1)}$,
which shows boundedness of the operator~$T^\delta$. 
In a similar manner, one can verify that $\|y^\delta\|_{H^{-1}(\Omega)} \le C'$.

As a next step, we verify the validity of the two conditions in~\eqref{eq:ass-noise}.  
By estimating the dual norm as before, and using the triangle inequality, we see that
\begin{align*}
\|T^\delta c - Tc \|_{H^{-1}(\Omega)} 
&= \|\div( c(\phi^\delta ) \nabla \phi^\delta ) - \div( c(\phi ) \nabla \phi )\|_{H^{-1}(\Omega)} \\
&\le \| c(\phi^\delta)\nabla \phi^\delta - c(\phi)\nabla \phi \|_{L^2(\Omega)} \\
&\le \| c(\phi^\delta) \nabla ( \phi^\delta - \phi) \|_{L^2(\Omega)} + \| ( c(\phi) -c(\phi^\delta)  )\nabla \phi \|_{L^2(\Omega)}.
\end{align*}
By the mean value theorem, we further obtain
\begin{align*}
\|T^\delta c - T c \|_{H^{-1}(\Omega)}
&\leq \| c \|_{L^\infty(-1,1)} \|\nabla (\phi^\delta - \phi) \|_{L^2(\Omega)} \\
&\qquad + \| c' \|_{L^\infty(-1,1)} \| \phi- \phi^\delta \|_{L^\infty(\Omega)} \|\nabla \phi\|_{L^2(\Omega)} 
\leq C \|\phi - \phi^\delta\|_{H^1(\Omega)}.
\end{align*}
In the last step, we used Sobolev embeddings and that fact that $c=c^\dag \in H^2(-1,1)$ is uniformly bounded by assumptions (A2)--(A3).
From \eqref{eq:ass-noise}, we thus conclude that 
\begin{align*} %\label{eq:estpertc}
    \|T^\delta c^\dag - Tc^\dag \|_{H^{-1}(\Omega)} 
    &\leq C \delta.
\end{align*}
In a similar manner, we can estimate the perturbations in the data by
\begin{align*}
\| y^\delta - y\|_{H^{-1}(\Omega)} 
\leq  \| \gamma \div(b(\phi^\delta)\nabla \Delta \phi^\delta ) &- \gamma \div(b(\phi)\nabla \Delta \phi ) \|_{H^{-1}(\Omega)} \notag \\
&  + \| \dt \phi^\delta - \dt \phi \|_{H^{-1}(\Omega)} 
\leq C' \delta. %\label{eq:estpertb}
\end{align*}
The results of Lemma~\ref{lem:convergenceapproximations} then guarantee stability and convergence of the regularized approximations $c_\alpha^\delta$ for the minimum norm solution $c^\dag$. 

\begin{remark}
If we have observations on a whole time interval $(t_1,t_2)$ it suffices to assume $\phi \in L^\infty (H^3(\Omega))$, where $L^p(X) = L^p(t_1,t_2;X)$ in the following,
and to require that
\begin{align} \label{eq:ass-noise2}
        \| \phi - \phi^\delta \|_{L^2(H^3(\Omega))} \leq \delta, \qquad \| \dt \phi - \dt \phi^\delta \|_{L^2(H^{-1}(\Omega))} \leq \delta,
\end{align}
which is slightly weaker than \eqref{eq:ass-noise}. 
By interpolation, we also obtain $\| \phi - \phi^\delta \|_{L^\infty(H^1(\Omega))} \leq c \delta$.
We then use $Y=L^2(H^{-1}(\Omega))$ as the image space of the operator $T^\delta$. 
By similar arguments as above, one can show that $\|T^\delta c\|_{L^2(H^{-1}(\Omega))}$, as well as
\begin{align*}
\|T^\delta c - Tc \|_{L^2(H^{-1}(\Omega))} \leq C \delta 
\qquad \text{and} \qquad 
\| y^\delta - y\|_{L^2(H^{-1}(\Omega))} 
\leq C' \delta;
\end{align*}
the details are left to the reader. Convergence of the regularized solutions $c_\alpha^\delta$ to a minimum-norm solution of the problem is again guaranteed by Lemma~\ref{lem:convergenceapproximations}.
\end{remark}

\subsection{Identification of $b(\cdot)$}
We now assume $\gamma$ and $f(\cdot)$ to be known and consider the identification of the mobility function $b(\cdot)$, for which we rewrite \eqref{eq:ch1} as
\begin{align} \label{eq:invb}
    \div\left( b(\phi)\nabla \mu \right) = \dt \phi.
\end{align}
From the data $\phi^\delta=\phi^\delta(\cdot,t)$ at $t \in (0,T)$, we can determine an approximation 
\begin{align*}
\mu^\delta = - \gamma \Delta \phi^\delta + f(\phi^\delta) 
\end{align*}
for the chemical potential mimicking equation~\eqref{eq:ch2}. 
Simply replacing $\phi$ and $\mu$ in equation~\eqref{eq:invb} by $\phi^\delta$ and $\mu^\delta$ then leads to a linear operator equation of the form \eqref{eq:linop} with 
perturbed operator 
\begin{align} \label{eq:Tb}
T^\delta: H^2(-1,1) \to H^{-1}(\Omega), \qquad b(\cdot) \mapsto \div (b( \phi^\delta) \nabla \mu^\delta ).
\end{align}
As before, we tacitly assumed $\phi^\delta(x,t) \in (-1,1)$.
With similar arguments as above, one can show that $T^\delta$ is linear and bounded, and that $y^\delta = \dt \phi^\delta$ lies in the data space $Y=H^{-1}(\Omega)$.
In order to verify \eqref{eq:condition}, we first estimate the perturbation in the chemical potential. By triangle inequality  and mean value theorem, we get
\begin{align*}
\|\mu^\delta - \mu \|_{H^1(\Omega)} 
&\leq \| -\gamma \Delta (\phi^\delta - \phi) +  f(\phi^\delta) - f(\phi)  \|_{H^1(\Omega)} \\
&\leq \gamma \|\phi^\delta - \phi \|_{H^3(\Omega)} + \| f'  \|_{L^\infty(-1,1)} \|\phi^\delta - \phi \|_{H^1(\Omega)} 
\le  c \delta.
\end{align*}
For the last step, we used uniform bounds for the parameters and assumption~\eqref{eq:ass-noise} on the data noise. 
Proceeding with similar arguments, we can then estimate
\begin{align*}
\|T^\delta b  - Tb \|_{H^{-1}(\Omega)} 
&\le \|b(\phi^\delta) \nabla \mu^\delta - b(\phi) \nabla \mu\|_{H^{-1}(\Omega)}\\
&\leq 
    \| b  \|_{L^\infty(-1,1)} \|\nabla ( \mu^\delta - \mu) \|_{L^2(\Omega)} \\
& \qquad \qquad + \| b' \|_{L^\infty(-1,1)} \|\nabla \mu \|_{L^2(\Omega)} \| \phi- \phi^\delta \|_{L^\infty(\Omega)}  
\leq C \, \delta.
\end{align*}
For the last step, we used Sobolev embeddings and assumed uniform bounds for the parameters $b \in H^2(-1,1)$ under consideration. 
From assumption~\eqref{eq:ass-noise}, we further conclude that $\|y - y^\delta\|_{H^{-1}(\Omega)} = \|\dt \phi - \dt \phi^\delta\|_{H^{-1}(\Omega)} \le \delta$. 
The results of Lemma~\ref{lem:convergenceapproximations} thus guarantee convergence of the regularized approximations $b_\alpha^\delta$ to a minimum norm solution $b^\dag$.

\begin{remark}
If we have measurements $\phi^\delta$ for a whole time interval $(t_1,t_2)$, we simply use $Y=L^2(H^{-1}(\Omega))$ as data space in the definition of the operator $T^\delta$, and assume the bounds \eqref{eq:ass-noise2} for the data noise. 
With the same arguments as above, one can then verify that $T^\delta$ is linear and bounded, as well as
\begin{align*}
\|T^\delta b  - Tb \|_{L^2(H^{-1}(\Omega))} \leq C  \, \delta 
\qquad \text{and} \qquad 
\|y - y^\delta\|_{L^2(H^{-1}(\Omega))} \le \delta.
\end{align*}
As a consequence, the abstract results of Lemma~\ref{lem:convergenceapproximations} can again be applied.
\end{remark}

\subsection{Simultaneous identification of $f(\cdot)$ and $b(\cdot)$} 
We assume $\gamma>0$ to be known and study the identification of $b(\cdot)$ and $f(\cdot)$. 
We start with rewriting the system \eqref{eq:ch1}--\eqref{eq:ch2} into 
\begin{align} \label{eq:invbc}
    - \gamma \div( b(\phi) \nabla \Delta \phi ) + \div( c(\phi)\nabla \phi ) = \dt \phi,
\end{align}
by eliminating $\mu$ using \eqref{eq:ch1} and introducing $c(s) = f'(s) b(s)$, as before.
Note that the function $f'(\cdot)$ can again be recovered from knowledge of $b(\cdot)$ and $c(\cdot)$.

Inserting data $\phi^\delta=\phi^\delta(\cdot,t)$ for a single time instance into this equation leads to a linear operator equation of the form \eqref{eq:linop}, with parameter $x=(b,c)$, data $y^\delta=\dt \phi^\delta$, and perturbed forward operator 
\begin{align} \label{eq:Tbc}
    T^\delta: H^2(-1,1)^2 &\mapsto H^{-1}(\Omega),\\
    (b,c) &\mapsto -\gamma \div( b(\phi^\delta)\nabla \Delta \phi^\delta ) + \div (c(\phi^\delta)\nabla \phi^\delta ). \notag
\end{align}
With the same arguments as used in the previous subsections, one can verify that $T^\delta$ is linear and bounded, as well as \begin{align*}
    \|T^\delta (b,c) - T(b,c) \|_{H^{-1}(\Omega)} \leq C\, \delta.
\end{align*}
Furthermore $\|y - y^\delta\|_{H^{-1}(\Omega)} \le \delta$ by assumption~\eqref{eq:ass-noise}. 
Therefore, convergence of the approximate solutions $(b_\alpha^\delta,c_\alpha^\delta)$ to a minimum-norm solution $(b,c)$ of the parameter identification problem is again ensured by Lemma~\ref{lem:convergenceapproximations}. 

\begin{remark}
If data $\phi^\delta$ are available on a whole time interval $[t_1,t_2]$, we can simply use $Y=L^2(H^{-1}(\Omega))$ as data space in the definition of the operator $T^\delta$ and utilize \eqref{eq:ass-noise2} as assumption on the data noise. 
This allows to establish the conditions \eqref{eq:condition} of Lemma~\ref{lem:convergenceapproximations} and yields a theoretical backup of the proposed regularization strategy.
\end{remark}

\section{Numerical illustration} 
\label{sec:num}

For illustration of our theoretical results, we briefly report on the actual performance of the proposed regularization strategies for a simple model problem. 
For sake of reproducability, we discuss in detail a one-dimensional test case, but note that similar results are obtained in two- and three space dimensions. 
Since more data are available in higher dimensions, the parameter reconstructions become more stable in higher dimensions. 

\subsection{Forward problem}
\label{sec:forward}
Let us start with describing the setup of ou rmodel problem:
As computational domain, we use $\Omega = (0,1) \approx \TT^1$, which is identified with the $1$-torus; hence \eqref{eq:ch1}--\eqref{eq:ch2} is actually supplemented by periodic boundary conditions.
We consider a polynomial double well potential 
\begin{align*}
    F(s) = (s-0.99)^2 (s+0.99)^4,
\end{align*}
and recall that only its derivative $f(s) = F'(s)$ appears in equation~\eqref{eq:ch2}.
As the mobility function for our model problem, we choose 
\begin{align*}
    b(\phi) = (1-\phi)^4 (1+\phi)^2 + 0.2,
\end{align*}
and we set $\gamma =0.003$ for the interface parameter.
Finally, the initial value for the phase fraction is prescribed by 
\begin{align*}
    \phi_0(x) = 0.1 \sin(2 \pi x  ) - 0.1 \sin(4 \pi x  ) +0.1 \sin(12 \pi x  ) +0.1.
\end{align*}
From Lemma~\ref{lem:weak}, we then deduce that $\phi$ is uniformly bounded on $\Omega \times [0,T]$, so that the functions $F(\cdot)$ and $b(\cdot)$ could be modified outsides this interval. Up to such modification, which do not affect our analysis, the assumption (A0)--(A3) are thus satisfied.

\subsection{Data generation}

In order to produce appropriate data for the inverse problem, we compute an approximate solution $\phi_{h,\tau}$ using the structure-preserving variational discretization method described in \cite{brunk2021relative}. 
This method is based on quadratic finite elements in space and a Petrov-Galerkin time-discretization with piecewise linear ansatz functions. 
For our numerical tests, we use uniform grids in space and time with mesh size $h = 5 \cdot 10^{-3}$ and time step $\tau = 2 \cdot 10^{-5}$. All simulations are performed up to time $T = 0.02$. 

To avoid inverse crimes, we use a different discretization strategy for the inverse problem, which is based on cubic splines in space and piecewise linear approximation in time, and we use spatial and temporal grids with the doubled mesh sizes. 
In a first step, we compute a corresponding approximation $\tilde\phi_{2h,2\tau}$ by interpolation of the data $\phi_{h,\tau}$. 
These data play the role of the perturbed data $\phi^\delta$ in our theoretical results; note that the perturbations here stem from discretization and interpolation errors. 
We further compute a cubic spline approximation $\tilde \mu_{2h,2\tau}$ by mimicking the identity \eqref{eq:ch2}.
\begin{figure}[ht!]
\centering
\footnotesize
\begin{tabular}{cc}
% \hspace{-0.7cm}
$\tilde \phi_{2h,2\tau}(x,t)$ & $  \partial_x \tilde \mu_{2h,2\tau}(x,t) $  \\
\includegraphics[trim={0.5cm 0.5cm 0.5cm 0.1cm},clip,scale=0.4]{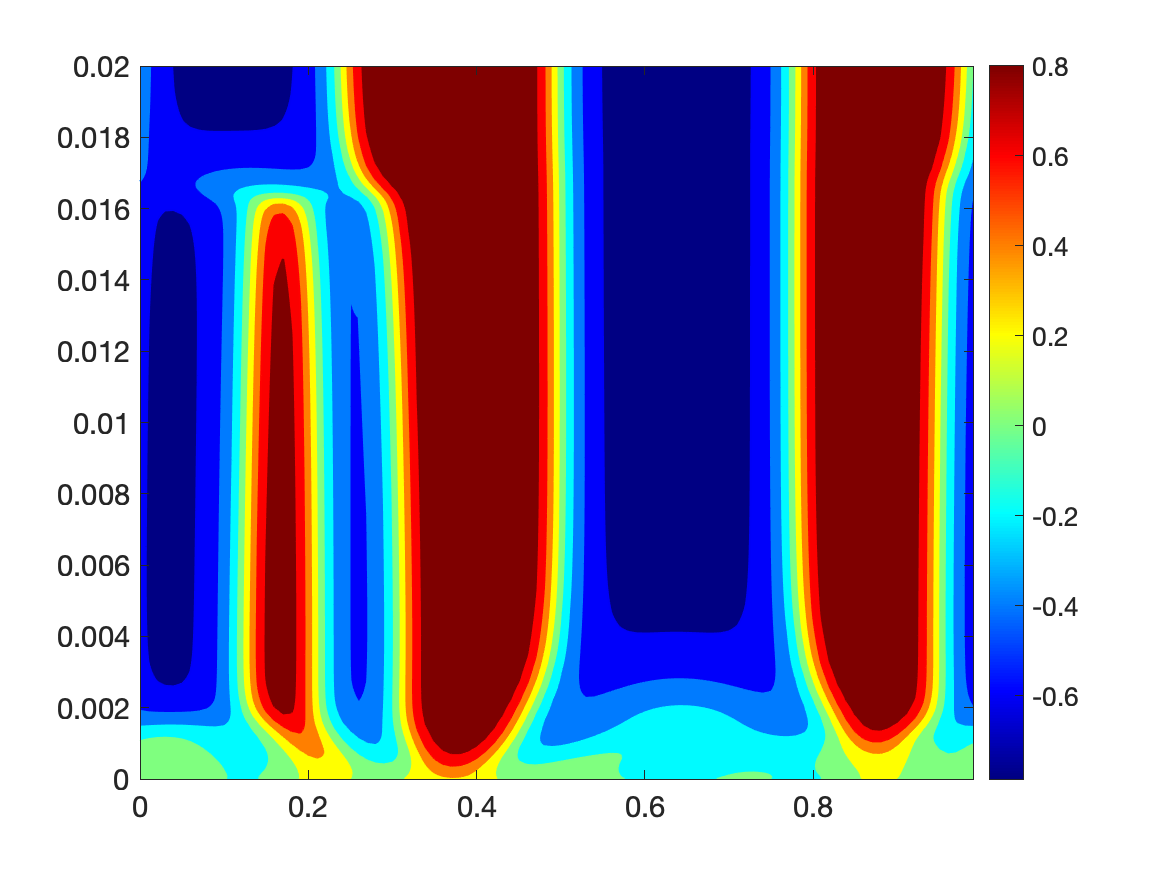} 
&
\includegraphics[trim={0.5cm 0.5cm 0.5cm 0.1cm},clip,scale=0.4]{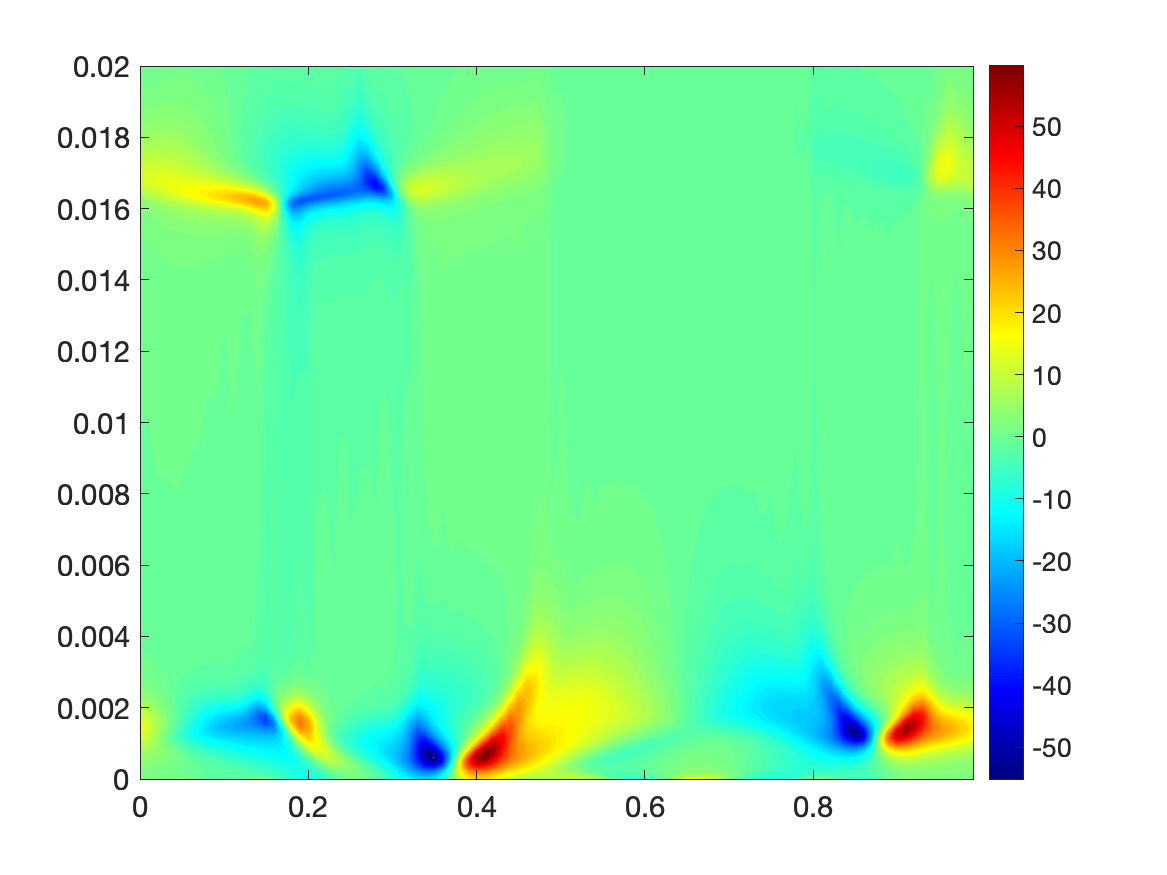}  
\end{tabular}
\caption{Contour plots of $\tilde \phi_{2h,2\tau}(x,t)$ (left) and $ \partial_x \tilde \mu_{2h,2\tau}(x,t)$ (right) with $x \in (0,1)$ on the x-axis and $t\in [0,0.02]$ on the y-axis.   \label{fig:identcond}
The color bar for the left plot shows the range $R = \{\phi^\delta(x,t) : x \in \Omega, \ t \in (0,T)\}$ of data that are actually attained. 
The right plot reveals areas, where information about the mobility function $b(\cdot)$ can be infered from the data.
}
\end{figure} 
In Figure \ref{fig:identcond}, we depict contour plots of the functions $\tilde \phi_{2h,2\tau}$ and $\dx \tilde \mu_{2h,2\tau}$. 
From these plots, we can infer information about the intervals $R_t=\{s=\phi^\delta(x,t): x \in \Omega\}$ and $\tilde R_t = \{s =  \phi(x,t): x \in \Omega, \ \dx \mu(x,t) \ne 0 \}$, 
where the parameter functions can be uniquely determined; see Section~\ref{sec:pid} for details.

\subsection{Numerical solution of the inverse problems}

For the implementation of the equation error methods introduced in the previous section, we use the following discretization strategy: 
The functions $\phi^\delta$, $\mu^\delta$ used in the definition of the operators $T^\delta$ in the previous section are approximated by the cubic splines $\tilde \phi_{2h,2\tau}$ and $\tilde \mu_{2h,2\tau}$, respectively. 
Backward difference quotients are used to approximate time derivatives.
The parameter functions $f(\cdot)$ and $b(\cdot)$ are discretized by natural cubic splines on a uniform grid of the interval $[-1,1]$ with grid size $\sigma=0.1$.
For clarity, we briefly also discuss in more detail the implementation of the perturbed version of~\eqref{eq:invf1}, which is used for the identification of the $f(\cdot)$. %
The  other two inverse problems are discretized in a similar manner. 

\subsection*{Numerical realization of identifying $f(\cdot)$}

The right-hand side $y^\delta$ of the inverse problem are approximated by a vector $\texttt{y}$, whose $i$th entry is computed by
\begin{align*}
\texttt{y}_i 
= (d_\tau \tilde \phi_{2h}, \tilde \psi_i)_{L^2(\Omega)} - \gamma (b(\tilde \phi_{2h}) \nabla \Delta \tilde \phi_{2h},\nabla \tilde \psi_i)_{L^2(\Omega)}
\end{align*}
where $\tilde \psi_i$ is the $i$th periodic cubic spline basis function, $\tilde \phi_{2h} = \tilde \phi_{2h,\tau}(\cdot,t)$ is the evaluation of the data at time $t$, and $d_\tau \tilde \phi_{2h} = \frac{1}{2\tau} (\tilde \phi_{2h,2\tau}(\cdot,t) - \tilde \phi_{2h,2\tau}(\cdot,t-2\tau))$ is the approximation for the time derivative by the backward difference quotient.
The matrix representation of the operator $T^\delta: c \mapsto \div(c(\tilde \phi_{2h}) \nabla \tilde \phi_{2h})$ is assembled by
\begin{align*}
    \texttt{T}_{ij} = -(\theta_j(\tilde \phi_{2h}) \nabla \tilde \phi_{2h}, \tilde \psi_i)_{L^2(\Omega)}
\end{align*}
where $\theta_j$ is the $j$th natural cubic spline basis function for the parameter $c(s) = \sum_j c_j \theta_j(s)$. 
We further define matrices $\texttt{M}_{ij} = (\tilde \psi_j, \tilde \psi_i)_{H^1(\Omega)}$ and $\texttt{R}_{ij}=(\theta_j,\theta_i)_{H^2(-1,1)}$, representing the scalar products on $H^1(\Omega)$ and $H^2(-1,1)$, respectively. 

The discretization of the Tikhonov functional \eqref{eq:leastsquaresprob} for problem \eqref{eq:invf1} is then given by 
\begin{align*}
     (\texttt{T} \texttt{c} - \texttt{y})^\top \texttt{M}^{-1} (\texttt{T} \texttt{c} - \texttt{y}) + \alpha \texttt{c}^\top \texttt{R} \texttt{c}. 
\end{align*}
Here $\texttt{c}=(c_1,\ldots,c_N)^\top$ is the coefficient vector of the parameter function $c$ to be determined. 
Minimization of this functional, respectively, solution of the corresponding normal equations can be achieved efficiently by the conjugate gradient algorithm.

\subsection{Numerical results}

We now briefly report on the results obtained for the three parameter identification problems discussed in Section~\ref{sec:pid} by the regularized equation error methods proposed in Section~\ref{sec:reg}. 
In all computations, the regularization parameter $\alpha$ is chosen heuristically, based on the L-curve-criterium; see e.g. \cite{engl1996book,hansen1993}.

\subsection*{Identification of $f(\cdot)$}

We assume $\gamma$ and $b(\cdot)$, defined as in Section~\ref{sec:forward}, to be known, and consider the identification of $f(\cdot)$. 
As outlined in Lemma~\ref{lem:invariance}, we can identify $f(\cdot)$ only up to a constant; hence only the derivative $f'(\cdot)$ can actually be identified.
Furthermore, the function $f'(\cdot)$ can obviously only be determined uniquely on the range
%$R_t=\{s=\phi(x,t) : x \in \Omega\}$ 
of data that are actually attained. 
\begin{figure}[ht!]
\centering
\footnotesize
\begin{tabular}{ccc}
% \hspace{-0.7cm}
$t=0.001$  & $t=0.008$ & $t \in [0,0.008]$ \\
\includegraphics[trim={1.5cm 0.9cm 1.5cm 0.6cm},clip,scale=0.28]{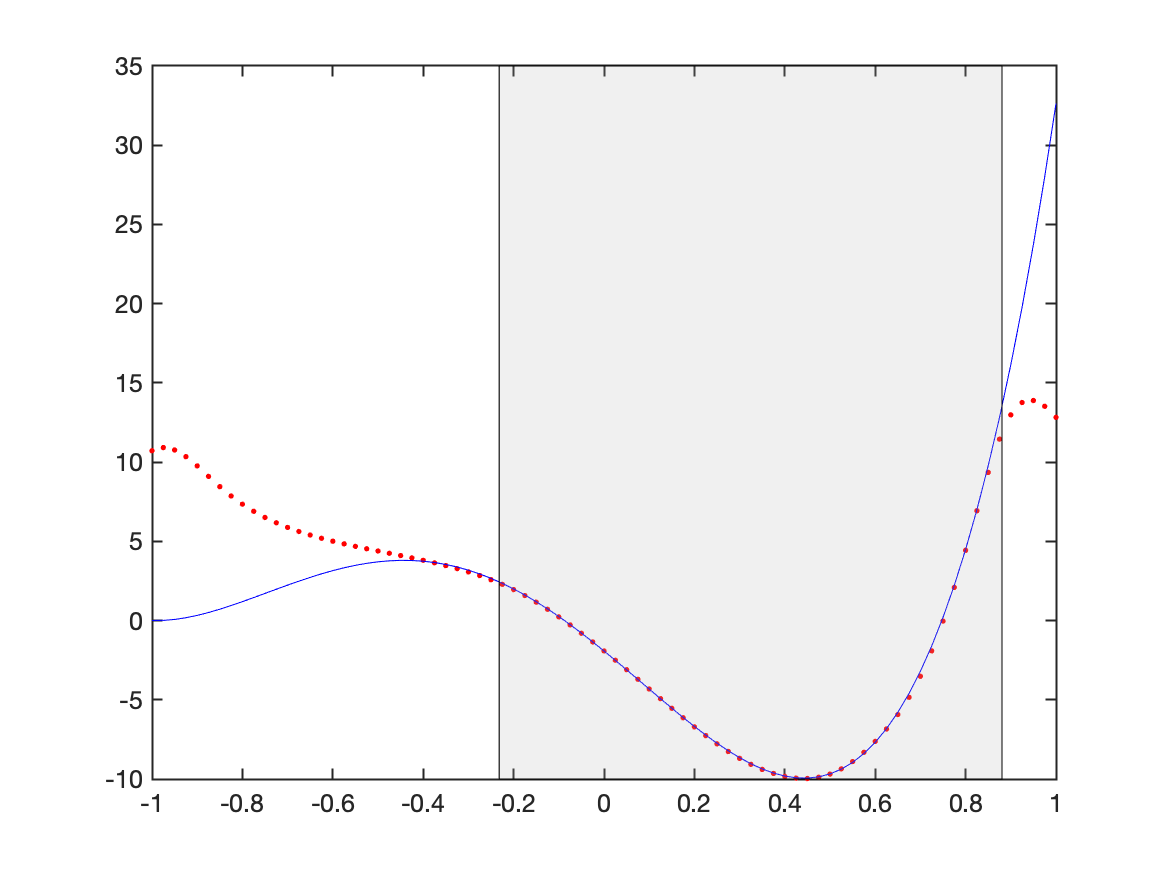} 
&
\includegraphics[trim={1.5cm 0.9cm 1.5cm 0.6cm},clip,scale=0.28]{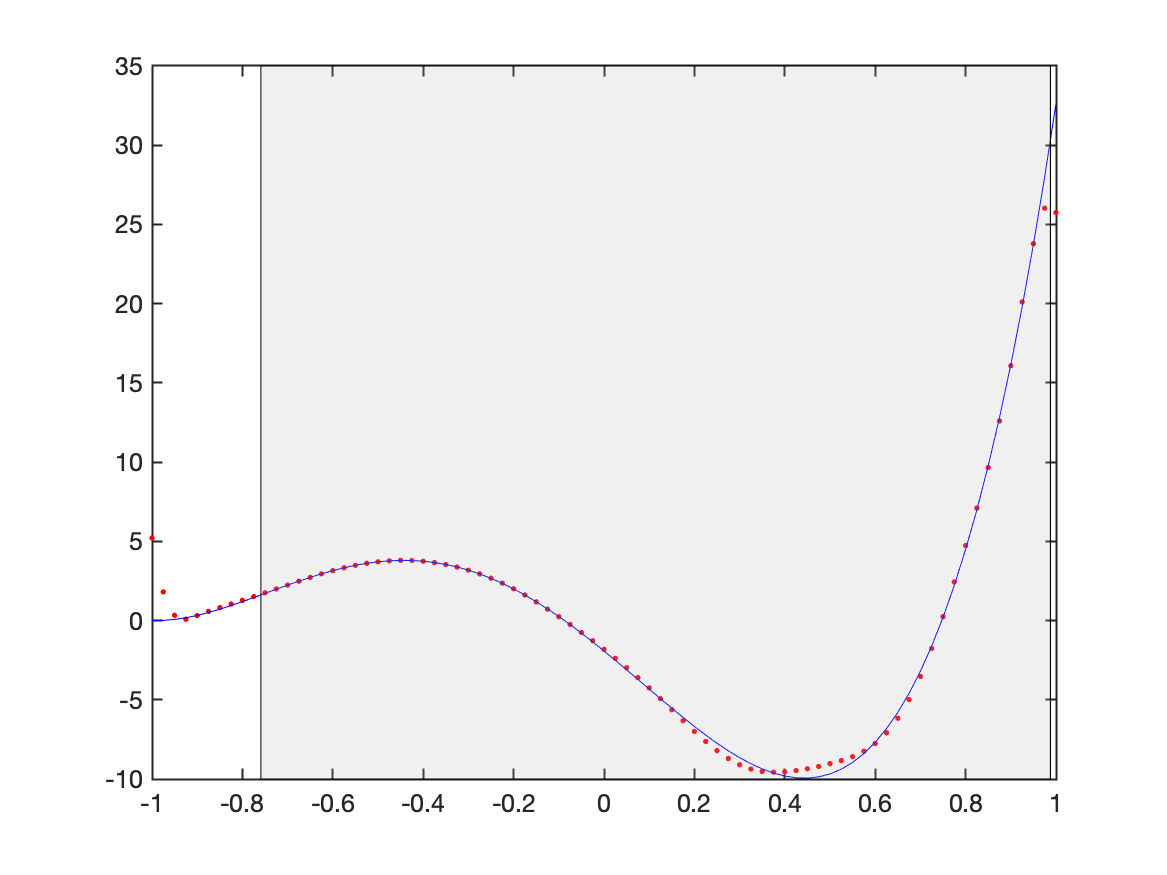}  
&
\includegraphics[trim={1.5cm 0.9cm 1.5cm 0.6cm},clip,scale=0.28]{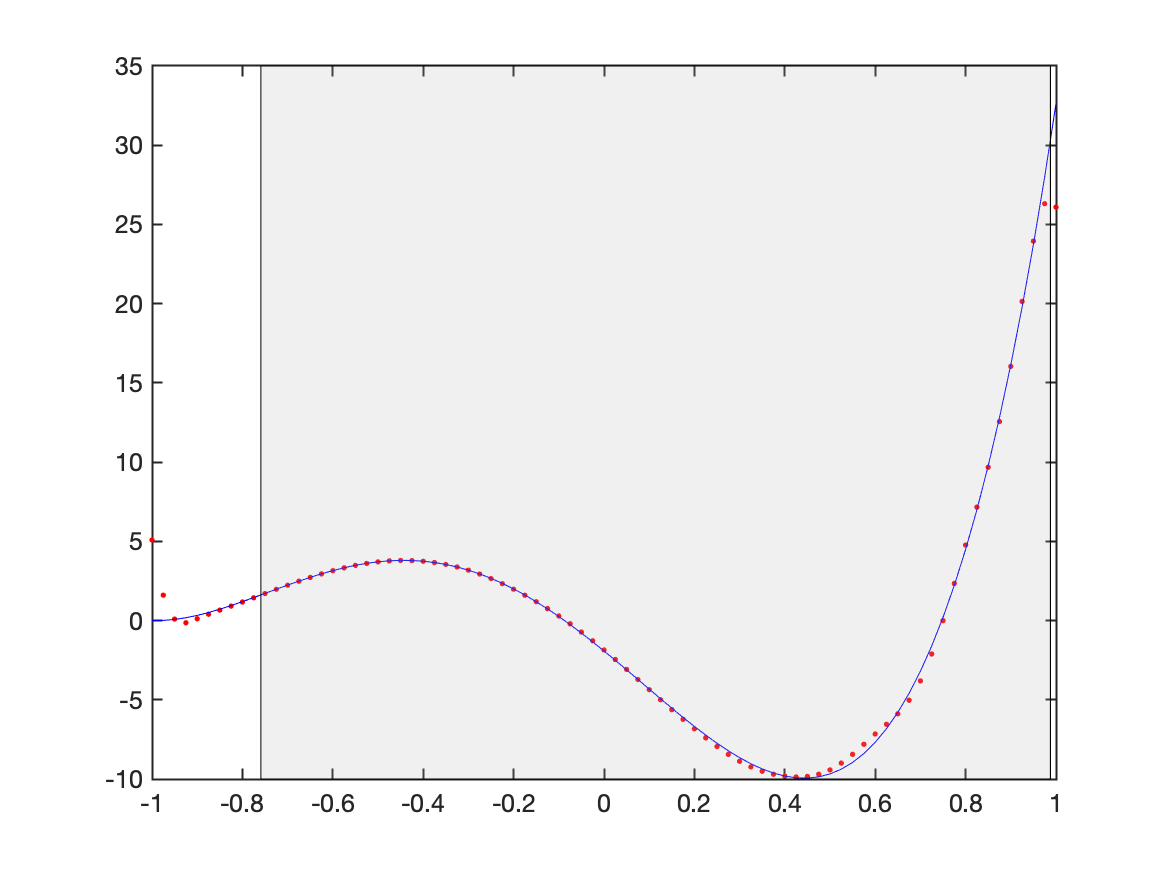} 
\end{tabular}
\caption{Reconstruction of $f'(\cdot)$ from data $\tilde \phi_{2h,2\tau}(\cdot,t)$, for $t$ as specified in the title of the plots.
The range $R_t$ of attained data is shaded in gray.
The solid blue line is the true function $f'(\cdot)$ while the dotted red line denotes the reconstruction $(f_\alpha^\delta)'(\cdot)$ obtained by the regularized equation error method of Section~\ref{sec:identify_f}. 
In all three test cases, the regularization parameter was chosen as $\alpha = 10^{-10}$. \label{fig:resultf}
}
\end{figure} 
In Figure~\ref{fig:resultf}, we display, the true value $f'(\cdot)=F''(\cdot)$ of the second derivative of the potential, and the corresponding reconstruction $(f_\alpha^\delta)'(\cdot)$, determined in our computations.  
The equation error method produces stable and accurate reconstructions in all cases.
As expected from Theorem~\ref{thm:uniq-f}, 
the function $f'(\cdot)$ is reconstructed reliably only on the respective range of data, while the regularization enforces stability but also a certain bias in the regions, where no data are available.

\subsection*{Identification of $b(\cdot)$}

We assume $f(\cdot)$ and $\gamma$ to be known, see Section~\ref{sec:forward}, and consider identification of the mobility. 
\begin{figure}[ht!]
\centering
\footnotesize
\begin{tabular}{ccc}% \hspace{-0.7cm}
$t=0.001$  & $t=0.002$ & $t \in [0,0.008]$ \\
\includegraphics[trim={1.5cm 0.9cm 1.5cm 0.6cm},clip,scale=0.28]{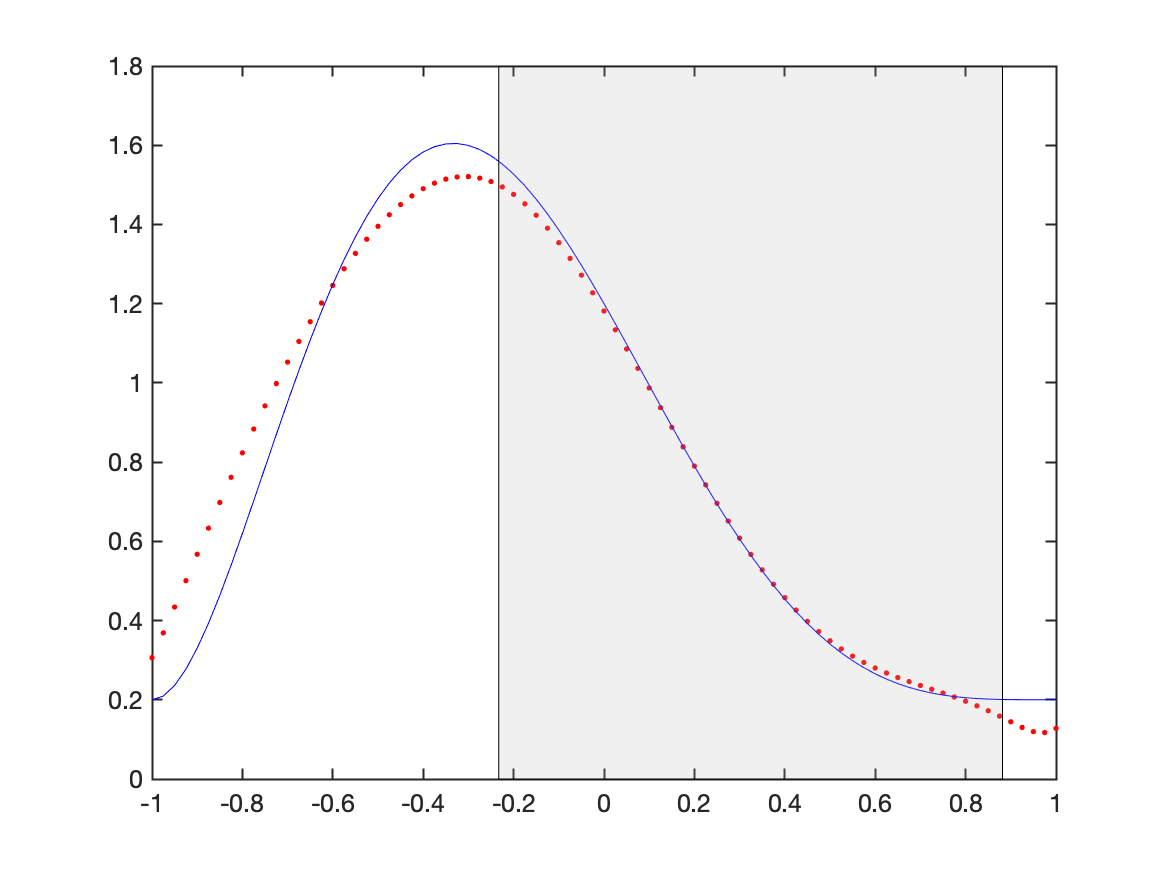} 
&
\includegraphics[trim={1.5cm 0.9cm 1.5cm 0.6cm},clip,scale=0.28]{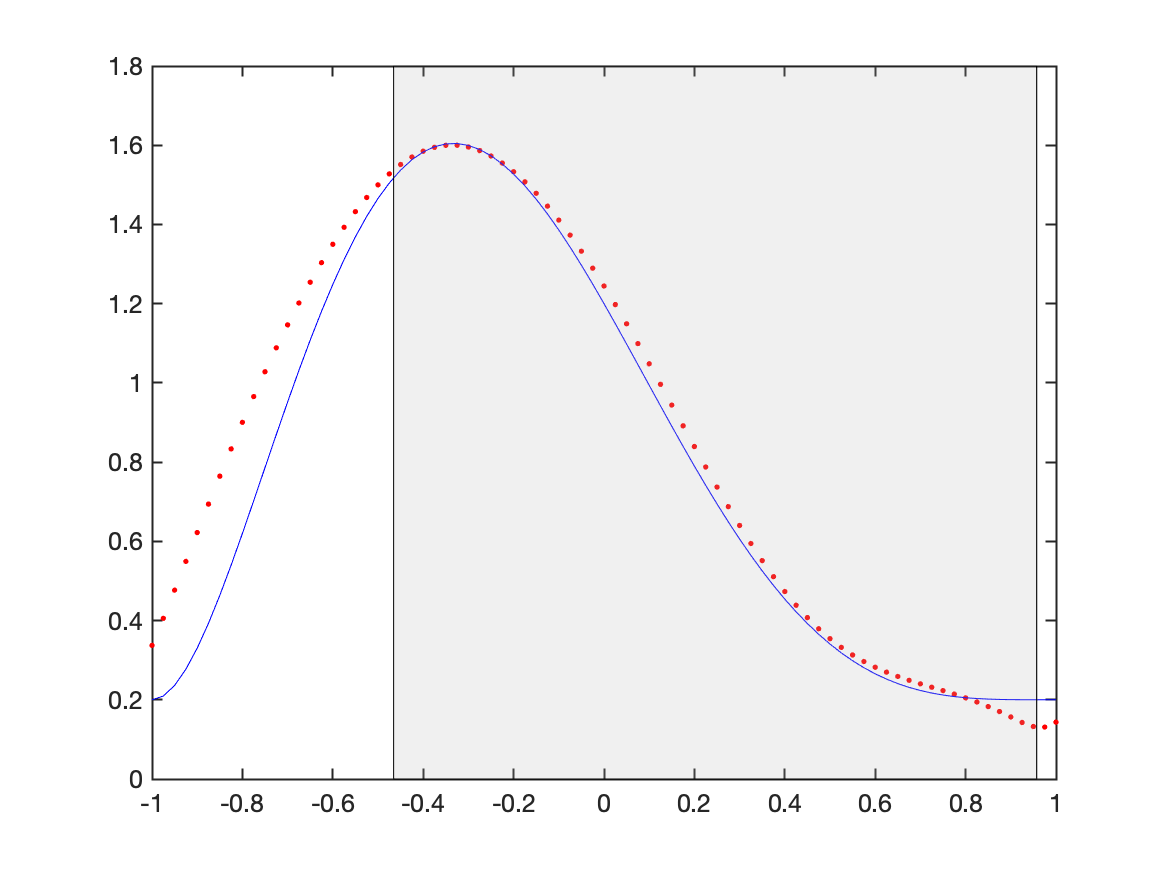}  
&
\includegraphics[trim={1.5cm 0.9cm 1.5cm 0.6cm},clip,scale=0.28]{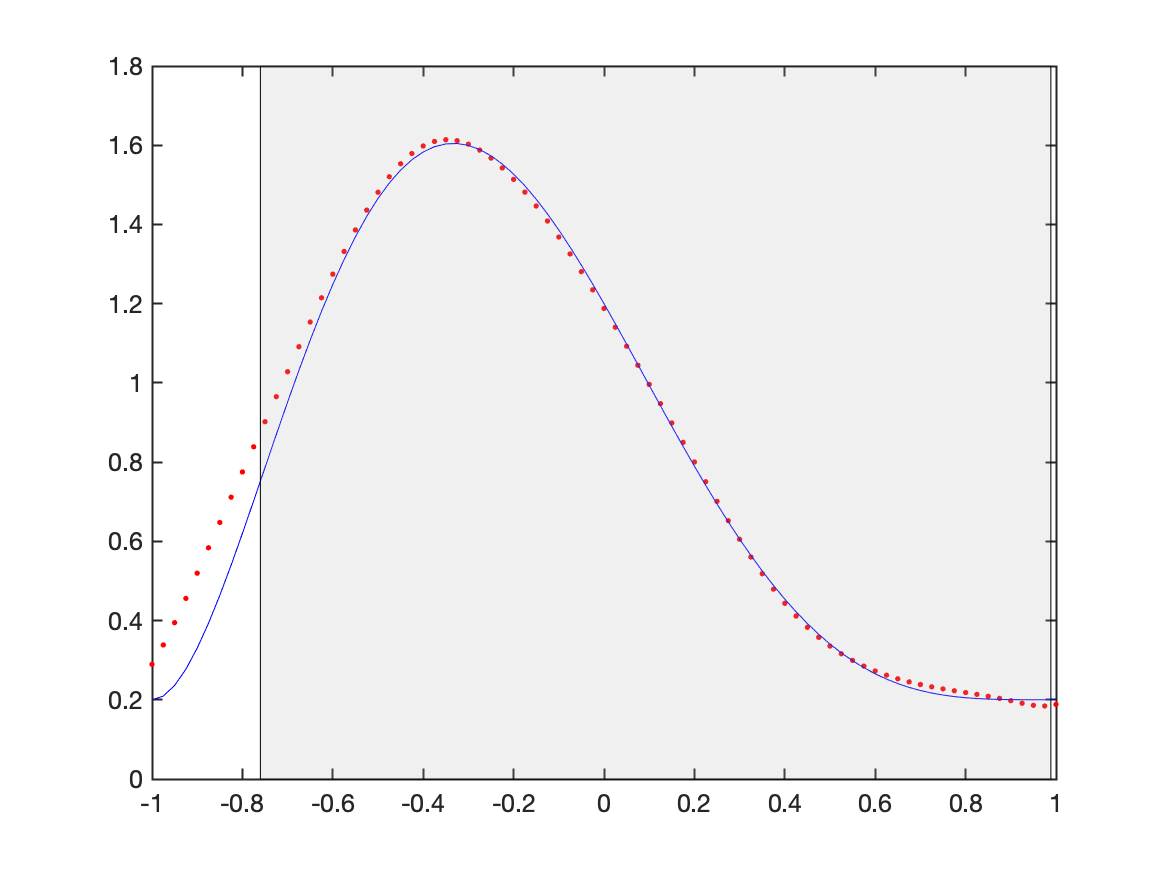} 
\end{tabular}
\caption{Reconstructions of the mobility function $b(\cdot)$ from perturbed data $\tilde \phi_{2h,2\tau}(\cdot,t)$ for time points $t$ depicted in the title of the plots.
The range $R_t$ of attained data is depicted in gray. 
The solid blue line depicts the true parameter, while the dotted red line denotes the reconstructions.  
The regularization parameters was set to $\alpha = 10^{-6}$ for all tests.\label{fig:resultb}
}
\end{figure} 
From Theorem~\ref{thm:identb}, we deduce that $b(\cdot)$ can be determined uniquely only on the range $\tilde R_t=\{s = \phi(x,t) : x \in \Omega, \ \dx \mu(x,t) \ne 0\}$ of the attained data, where the gradient of the chemical potential does not vanish. 
In Figure~\ref{fig:resultb}, we display the reconstructions obtained from distributed phase fraction data data for single time steps and for a whole time interval. 
Again, the reconstructed mobility is in good agreement with the true parameter $b(\cdot)$ on the range of attained data, while the reconstructions outside this range are stable but biased by the regularization term in the equation error method.

\subsection*{Simultaneous identification of $f(\cdot)$ and $b(\cdot)$}

Here we only assume $\gamma$ to be known and chosen as in Section~\ref{sec:forward}.
According to Theorem~\ref{thm:bc}, the simultaneous identification of both parameters requires data at multiple time steps. 
\begin{figure}[ht!]
\centering
\footnotesize
\begin{tabular}{cc}
% \hspace{-0.7cm}
$b(\cdot)$  & $f'(\cdot)$   \\
\includegraphics[trim={1.5cm 0.9cm 1.5cm 0.6cm},clip,scale=0.34]{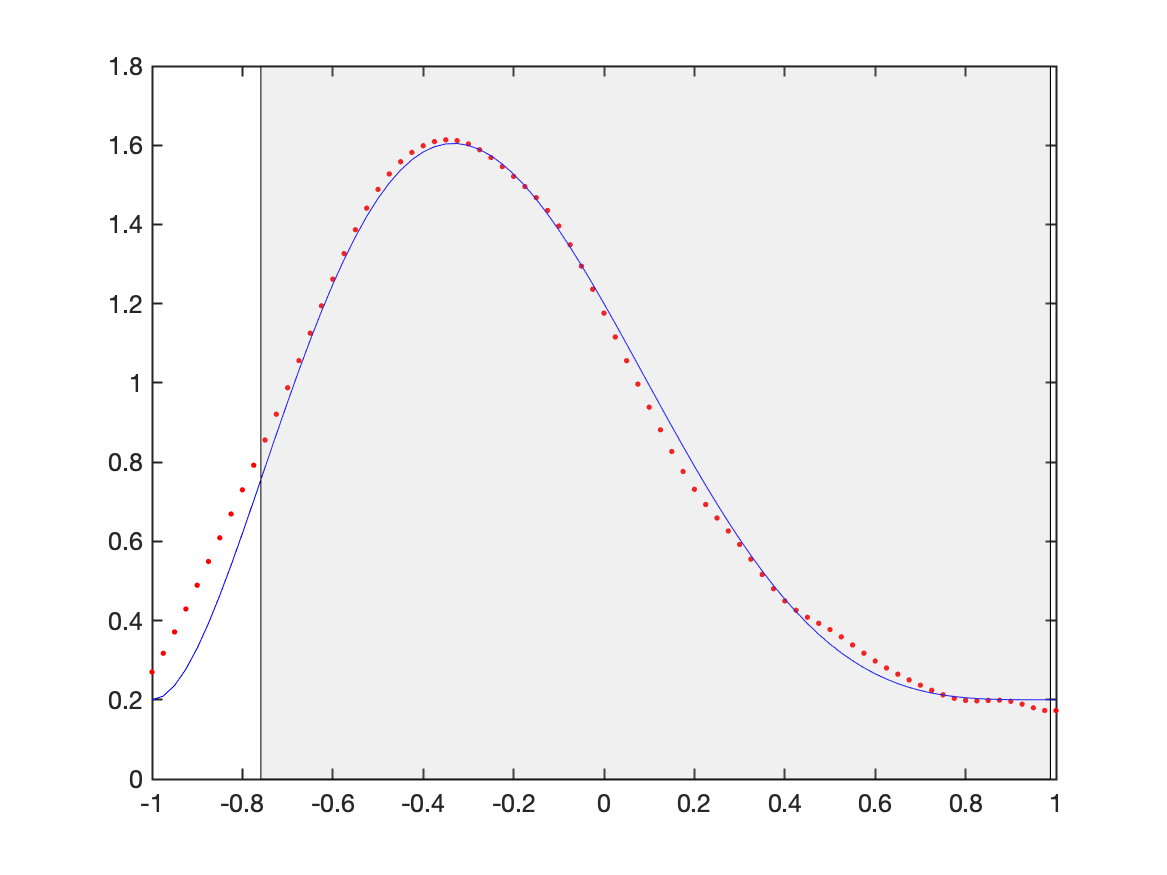} 
&
\includegraphics[trim={1.5cm 0.9cm 1.5cm 0.6cm},clip,scale=0.34]{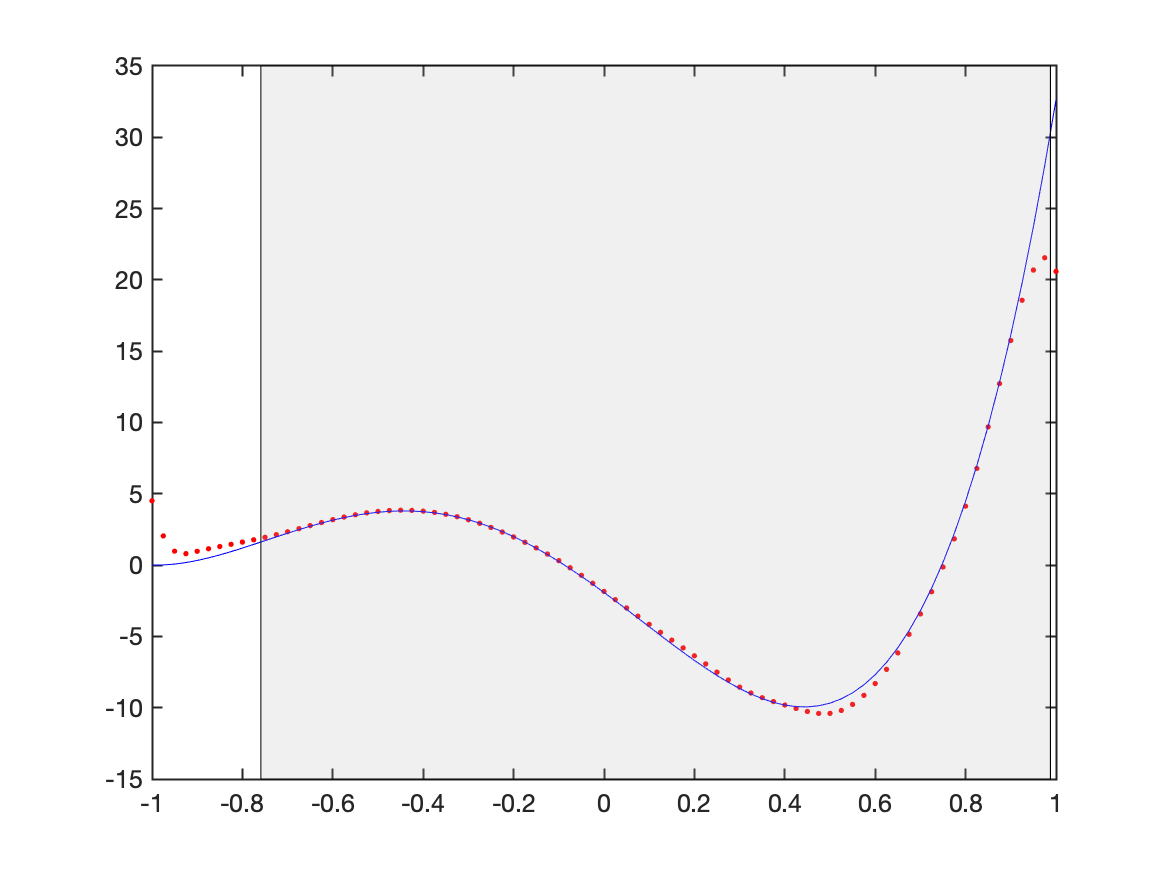} 
\end{tabular}
\caption{Simultaneous reconstructions of $b(\cdot)$ and $f'(\cdot)$ from perturbed data $\tilde \phi_{2h,2\tau}(\cdot,t)$ with $t \in [0,0.008]$.
The range of the attained data is again depicted in gray.
The solid blue line depicts the true parameter functions, while the corresponding reconstructions are denoted by dotted red lines.
The regularization parameter was chosen as $\alpha = 10^{-9}$ here. \label{fig:resultfb}
}
\end{figure} 
In Figure~\ref{fig:resultfb}, we therefore only report about reconstructions obtained for data on a whole time interval. 
As expected, the parameter functions are stably and accurately determined on the range of available data.
In our numerical tests, we also checked the validity of the observability condition stated in Theorem~\ref{thm:bc}, by numerical computing the parameters $A_b(s,t_i)$, $A_c(s,t_i)$ for different values of $s$ and $t_i$,
and observed the required linear independence in most cases. 

\subsection*{Multi-dimensional problems}

We also performed numerical tests for a similar model problem in $2d$ and observed very similar results to that for the $1d$ problem presented above. 
In fact, the reconstructions  in $2d$ are more stable, which is not surprising, since more data are available, while the functions to be determined, stay the same.

\section{Conclusion}

In this paper, we studied the identification of multiple parameter functions in a Cahn-Hilliard model for phase-separation from a theoretical and a numerical point of view.
Identifiability up to certain invariances was established under mild and natural observability conditions on the data which, in principle, can be verified prior to the computations.
A regularized equation error approach was studied for the numerical solution and viability of the parameter estimation in the presence of discretization errors was demonstrated. 

Various generalizations of the Cahn-Hilliard system have been proposed to derive more realistic models for phase separation processes; see \cite{Ipocoana,Marveggio} or \cite{Abels} for examples.
With increasing complexity, also more model parameters or parameter functions are introduced and have to be calibrated to obtain quantitative agreement with microscopic simulations or experimental data. The parameter identification for such problems may be substantially more involved and will be investigated in future research.

\section*{Acknowledgement}
Financial support of this work by the German Science Foundation (DFG) via grants TRR~146 (project~C3) and SPP~2256 (project Eg-331/2-1) is gratefully acknowledged.

%\bibliographystyle{abbrv}
%\bibliography{bibfile}

\end{document}